\documentclass[12pt]{article}
\usepackage{amsmath}
\usepackage{amsfonts}
\usepackage{amssymb}
\usepackage{amsthm}
\theoremstyle{definition}
\newtheorem{theorem}{Theorem}
\newtheorem{definition}{Definition}

\begin{document}
\title{Eigenstates of C*-Algebras}
\author{Luther Rinehart}
\maketitle
\begin{abstract}
	We introduce the notion of eigenstate of an operator in an abstract C*-algebra, and prove several properties. Most significantly, if the operator is self-adjoint, then every element of its spectrum has a corresponding eigenstate.
\end{abstract}
In quantum mechanics an important role is played by eigenvectors of operators, namely those vectors $v$ such that $x(v)=\lambda v$ for some $\lambda\in \mathbb{C}$. Such vectors represent physical states in which the operator $x$ has the definite value $\lambda$. The algebraic formulation of quantum mechanics seeks to develop the theory using only a $C^*$-algebra of operators, without reference to a Hilbert space. The purpose of this paper is to examine the analogous concept to eigenvector in the purely algebraic setting.

A key advantage of considering algebraic states is shown by theorem \ref{spectrum}. Typically, not every point in an operator's spectrum has a corresponding eigenvector. Indeed, an operator might have no eigenvectors at all. However, theorem \ref{spectrum} shows that, as long as the operator is self-adjoint, every point in the spectrum has an algebraic eigenstate.

In the following, $A$ will always be a $C^*$-algebra with 1. Its space of states $S(A)$ consists of the normalized positive linear functionals.
\begin{definition}
Let $x\in A$. A state $E\in S(A)$ is an \emph{eigenstate} of $x$ with \emph{eigenvalue} $\lambda\in\mathbb{C}$ if $\forall y\in A,\ E(yx)=\lambda E(y)$. 
\end{definition}
\begin{theorem}
	Every eigenvalue of $x$ is contained in the spectrum of $x$.
	\begin{proof}
		Let $\lambda$ be an eigenvalue of $x$. For contradiction, suppose $y=(x-\lambda)^{-1}$. Let $E$ be the associated eigenstate for $\lambda$. By definition, $E(y(x-\lambda))=0$. But this says $E(1)=0$, a contradiction. So $x-\lambda$ is not invertible.
	\end{proof}
\end{theorem}
\begin{theorem}\label{var}
A state $E$ is an eigenstate of $x$ with eigenvalue $\lambda$ if and only if $E((x-\lambda)^*(x-\lambda))=0$.
\begin{proof}
The implication $\Rightarrow$ follows by direct computation or by noting that since $\forall y,\ E(y(x-\lambda))=0$, then this certainly holds for $y=(x-\lambda)^*$.\\
\\
To show the converse, apply the Cauchy-Schwartz inequality for states:
$$\forall y\ |E(y(x-\lambda))|^2\leq E(y^*y)E((x-\lambda)^*(x-\lambda))=0$$
\end{proof}
\end{theorem}
The next theorem shows that the notion of eigenstate coincides with the standard notion of eigenvector if we fix a representation and focus only on vector states.
\begin{theorem}
Let $\pi\colon A\rightarrow B(H)$ be a *-representation of $A$ and let $v\in H$. Let $E_v$ be the vector state associated with $v$, namely, $E_v(x)=\langle v,\pi(x)v \rangle$. $E_v$ is an eigenstate of $x$ with eigenvalue $\lambda$ if and only if $\pi(x)v=\lambda v$.
\begin{proof}
Suppose $\pi(x)v=\lambda v$. $\forall y,\ E_v(yx)=\langle v, \pi(yx)v \rangle=\langle v, \pi(y)\pi(x)v \rangle=\lambda \langle v,\pi(y)v \rangle=\lambda E_v(y)$.\\
\\
Conversely, suppose $E_v$ is an eigenstate of $x$ with eigenvalue $\lambda$. By theorem \ref{var}, $0=E_v((x-\lambda)^*(x-\lambda))=\langle v,\pi((x-\lambda)^*(x-\lambda))v \rangle =\langle \pi(x-\lambda)v,\pi(x-\lambda)v \rangle$. So $\pi(x-\lambda)v=0$, that is $\pi(x)v=\lambda v$.
\end{proof}
\end{theorem}
For the following theorem, recall that if $x=x^*$, the C*-algebra generated by $x$ can be canonically identified with $C(\sigma(x))$, the continuous functions on the spectrum of $x$.
\begin{theorem}
Let $x=x^*$. If $E$ is an eigenstate of $x$ with eigenvalue $\lambda$, then $\forall f\in C(\sigma(x))$, $E$ is an eigenstate of $f$ with eigenvalue $f(\lambda)$.
\begin{proof}
By induction on $n$, the theorem holds for $f=x^n$. Indeed, $\forall y,\ E(yx^n)=E(yx^{n-1}x)=\lambda E(yx^{n-1})=\lambda\lambda^{n-1}E(y)=\lambda^n E(y)$.\\
\\
By linearity, for all polynomials $p$, $E$ is an eigenstate of $p$ with eigenvalue $p(\lambda)$. Now let $f\in C(\sigma(x))$. By the Stone-Weierstrass theorem, there exists a sequence of polynomials $\{p_n\}$ converging uniformly to $f$. Using continuity of multiplication and of $E$, $\forall y,\  E(yf)=\lim E(yp_n)=\lim p_n(\lambda) E(y)=f(\lambda)E(y)$.
\end{proof}
\end{theorem}
If $x$ is self-adjoint, every point of its spectrum has a corresponding eigenstate. The proof is analogous to that of theorem 1.3.1 in \cite{murphy}.
\begin{theorem}\label{spectrum}
Let $x=x^*$. If $\lambda\in\sigma(x)$, then $\lambda$ is an eigenvalue of x.
\begin{proof}
Let $x=x^*$ and let $\lambda\in\sigma(x)$. Observe that finding an eigenstate with eigenvalue $\lambda$ is equivalent to finding a state which vanishes on the left ideal generated by $(x-\lambda)$. Let $J=A(x-\lambda)$ be this ideal. Since $(x-\lambda)$ is self-adjoint and not invertible, it is not left-invertible, so $1\notin J$. Furthermore, no element of $J$ can be invertible, as this would imply $y\in J\Rightarrow y^{-1}y=1\in J$. Thus $$\forall y\in J,\ 0\in\sigma(y)$$
$$\Rightarrow 1\in\sigma(y-1)$$
$$\Rightarrow 1\leq\|y-1\|$$
So $1\notin \bar{J}$. By the Hahn-Banach theorem, we can find a linear functional $E\in A^*$ satisfying $E\rvert_{\bar{J}}=0$ and $E(1)=1=\|E\|$. Such $E$ will be an eigenstate of $x$ with eigenvalue $\lambda$. 
\end{proof}
\end{theorem}
\begin{theorem}
Any set of eigenstates of $x$ all having distinct eigenvalues, are linearly independent.
\begin{proof}
Let $\{E_\alpha\}$ be a collection of eigenstates of $x$ with eigenvalues $\{\lambda_\alpha\}$, and $\{\lambda_\alpha \}$ all distinct. By induction on $n$, if $\sum_{\alpha=1}^{n} c_\alpha E_\alpha =0$, then $c_1=...=c_n=0$. \\
\\
The case $n=1$ is clear.\\
\\
Make the induction hypothesis that the conclusion holds for $n=k$, and suppose that $\sum_{\alpha=1}^{k+1} c_\alpha E_\alpha =0$. Apply this sum to the element $y(x-\lambda_{k+1})$, with $y$ arbitrary:
$$0=\left( \sum_{\alpha=1}^{k+1}c_\alpha E_\alpha \right)  \left( y(x-\lambda_{k+1})\right) $$
$$=\sum_{\alpha=1}^{k}c_\alpha(\lambda_\alpha-\lambda_{k+1})E_\alpha(y)$$
By the induction hypothesis, $c_\alpha(\lambda_\alpha-\lambda_{k+1})=0$ for $\alpha=1...k$. But the $\{\lambda_\alpha\}$ are all distinct, so $\lambda_\alpha-\lambda_{k+1}\neq 0$ for $\alpha=1...k$. So $c_1=...=c_k=0$. Finally, $c_{k+1}=0$ by the $n=1$ case.
\end{proof}
\end{theorem}
\begin{definition}
A pair of positive linear functionals $\omega,\phi\in A^*_+$ are \emph{orthogonal} if $\|\omega-\phi\|=\|\omega\|+\|\phi\|$. (See definition 3.2.3 in \cite{pedersen})
\end{definition}
\begin{theorem}
If $x=x^*$, then any two eigenstates of $x$ with different eigenvalues are orthogonal. 
\begin{proof}
Let $E_1$ and $E_2$ be eigenstates of $x$ with respective eigenvalues $\lambda_1$ and $\lambda_2$, with $\lambda_1\neq \lambda_2$. Then there exists a continuous function $f\in C(\sigma(x))$ such that $\|f\|=1$, $f(\lambda_1)=1$, and $f(\lambda_2)=-1$. Then $(E_1-E_2)(f)=1-(-1)=2$. Consequently, $2\leq\|E_1-E_2\|\leq\|E_1\|+\|E_2\|=2$, so $E_1$ and $E_2$ are orthogonal.
\end{proof}
\end{theorem}
Finally, eigenstates of projection operators are precisely those obtained by projection.
\begin{theorem}
Let $p$ be a self-adjoint projection and $E$ a state. Let $pE$ be the state given by
$$pE(x)=\frac{E(pxp)}{E(p)}.$$
The following are equivalent:
\begin{enumerate}
\item $pE=E$
\item $E$ is an eigenstate of $p$ with eigenvalue 1
\item $E(p)=1$
\end{enumerate}
\begin{proof}
$(1)\Rightarrow (2)$: $\forall y\in A,\ pE(yp)=E(pypp)/E(p)=E(pyp)/E(p)=pE(y)$, so $pE$ is an eigenstate of $p$ with eigenvalue 1. \\
\\
$(2)\Rightarrow(3)$: Obvious.\\
\\
$(3)\Rightarrow(2)$: Let $E(p)=1$. $E((p-1)^*(p-1))=E((p-1)^2)=E(1-p)=1-E(p)=0$ and apply theorem \ref{var}.\\
\\
$(2)\Rightarrow(1)$: Suppose $E$ is an eigenstate of $p$ with eigenvalue 1, that is $\forall y\in A,\ E(yp)=E(y)$. By taking adjoint, this also gives $E(py^*)=E(y^*)$. Since $*$ is an involution, we have $\forall y\ E(py)=E(y)$. Now replacing $y$ with $yp$, $E(pyp)=E(yp)=E(y)$. And since $E(p)=1$, altogether this shows that $\forall y\ pE(y)=E(y)$, so $E=pE$.
\end{proof}
\end{theorem}

\end{document}